\documentclass[reqno]{amsart}

\usepackage{amssymb}
\usepackage{amsmath}
 \usepackage{amsthm}

\usepackage{amssymb,verbatim}
\usepackage{color}
\usepackage{url}
\usepackage{enumitem}
\usepackage{capt-of} 
\usepackage{hyperref}
\usepackage{graphicx} 
\hypersetup{colorlinks}
\usepackage{extarrows}
\usepackage[lined,boxed,norelsize,linesnumbered]{algorithm2e} 

\setlist[enumerate]{
label=\textnormal{({\roman*})},
ref={\roman*}}

\makeatletter 
\newtheorem*{rep@theorem}{\rep@title}
\newcommand{\newreptheorem}[2]{%
\newenvironment{rep#1}[1]{%
 \def\rep@title{#2 \ref{##1}}%
 \begin{rep@theorem}}%
 {\end{rep@theorem}}}
\makeatother

\newreptheorem{theorem}{Theorem} 
\newreptheorem{conjecture}{Conjecture} 
\newreptheorem{corollary}{Corollary}
\newreptheorem{lemma}{Lemma}
\newreptheorem{proposition}{Proposition}

\newtheorem{theorem}{Theorem}[section]
\newtheorem{lemma}[theorem]{Lemma}
\newtheorem{proposition}[theorem]{Proposition}

\newtheorem{corollary}[theorem]{Corollary}

\newtheorem{question}[theorem]{Question}

\theoremstyle{definition}
\newtheorem*{remarkx}{Remark}
\newenvironment{remark}
  {\pushQED{\qed}\remarkx}
  {\popQED\endremarkx}

\newtheorem{definitionx}[theorem]{Definition}
\newenvironment{definition}
  {\pushQED{\qed}\definitionx}
  {\popQED\enddefinitionx}

\DeclareMathOperator{\outdeg}{\textnormal{outdeg}}
\DeclareMathOperator{\indeg}{\textnormal{indeg}}
\newcommand{\Sand}{\textnormal{Sand}} 
\newcommand{\lvl}{\textnormal{lvl}} 
\newcommand{\cc}{\mathbf{c}}
\newcommand{\dd}{\mathbf{d}}
\newcommand{\ee}{\mathbf{e}}
\newcommand{\hh}{\mathbf{h}}
\newcommand{\cs}{\widehat{\mathbf{c}}}
\newcommand{\ds}{\widehat{\mathbf{d}}}
\newcommand{\qq}{\mathbf{q}}
\newcommand{\rr}{\mathbf{r}}
\newcommand{\us}{\widehat{\mathbf u}} 
\newcommand{\zz}{\mathbf{z}}
\newcommand{\satu}{\textbf{{1}}} 
\newcommand{\el}{{\sim}} 
\newcommand{\es}{\stackrel{\text s}{\sim}} 
\newcommand{\dar}{\longrightarrow} 

\newcommand{\N}{\mathbb{N}} 
\newcommand{\Z}{\mathbb{Z}} 
\newcommand{\R}{\mathbb{R}} 
\newcommand{\Rec}{\textnormal{Rec}} 

\newcommand{\LL}{\Delta} 
\newcommand{\Park}{\textnormal{Park}} 
\newcommand{\Rc}{\mathcal R} 
\newcommand{\B}{\mathcal B} 

\newcommand{\ext}{\textnormal{ext}} 
\newcommand{\BE}{\textnormal{BE}} 
\newcommand{\BV}{\textnormal{BV}} 

\newcommand{\src}{\textnormal{src}} 
\newcommand{\trgt}{\textnormal{trgt}} 

\begin{document}



\title[Abelian sandpile and Biggs-Merino polynomial for digraphs]{Abelian sandpile model  and Biggs-Merino polynomial for directed graphs}

 \author{Swee Hong Chan}
 \address{Department of Mathematics, Cornell University, Ithaca, NY 14853.}
\email{\url{sc2637@cornell.edu}}
\urladdr{https://www.math.cornell.edu/~sc2637/} 


\begin{abstract}
We prove several  results concerning a polynomial that arises from the sandpile model on directed graphs;  these results are previously only known for undirected graphs.
Implicit in the sandpile model is the choice of a sink vertex,
and it is conjectured by Perrot and Pham that the polynomial $c_0+c_1y+\ldots c_n y^n$, where $c_i$ is the number of recurrent classes of the sandpile model with level $i$, is independent of the choice of the sink.
 We prove their conjecture by expressing the polynomial as an invariant of the sinkless sandpile model.
 We then  present a bijection between  arborescences of directed graphs 
 and reverse $G$-parking functions that preserves external activity 
 by generalizing  Cori-Le Borgne bijection for undirected graphs.
As an application of this bijection, 
we extend Merino's Theorem by showing 
that for Eulerian directed graphs the polynomial $c_0+c_1y+\ldots c_n y^n$
is equal to the greedoid polynomial of the graph.
\end{abstract}



\maketitle


%
%
%
%
%
%
%
%

\section{Introduction}

To what extent do the known results for undirected graphs extend to directed graphs?
Driven by this question,
we consider a remarkable  theorem of Merino L{\'o}pez~\cite{Mer97} that expresses 
a one variable  specialization of the Tutte polynomial of an undirected 
graph in terms of the abelian sandpile model on the graph.
In this paper, we show that this theorem can be extended to all Eulerian directed graphs, and a weaker version of the theorem can be extended to all directed graphs.

The \emph{abelian sandpile model} is a dynamical system on a finite directed graph that  starts with a number of chips at each vertex of the graph.
If a vertex has 
 at least as many chips as its outgoing edges, then we are allowed to \emph{fire} the vertex  
by sending one chip along each edge leaving the vertex to the neighbors of the vertex. 
This  model was  introduced by Dhar \cite{Dhar90} as a model to study the concept of self-organized criticality introduced in \cite{BTW88}.  
Since then it has been studied in several different field of mathematics.
In graph theory it was studied under the name   of {chip-firing game}~\cite{Tar88, BLS91};
it appears in arithmetic geometry in the study of the Jacobian of algebraic curves \cite{Lor89, BN07};
 and in  algebraic graph theory it relates to the study of potential theory on graphs \cite{Biggs97pt, BS13}.
  

It is common to study  the  sandpile model   by specifying a vertex  as the sink vertex, and all chips that end up at the sink vertex are  removed from the process.
For a strongly connected directed graph, this guarantees that  
the sandpile model terminates (i.e. when 
  none of the vertices have enough chips to be fired) in finite time.
After fixing a sink,
one can study a special type of chip configurations with the following property.
A chip configuration 
 is \emph{recurrent} if, for any arbitrary chip configuration  as the initial state,
one can add a finite 
amount of chips to each vertex so that the recurrent configuration is the  state of the sandpile model when the process terminates.

It was conjectured by Biggs~\cite{Biggs97} and was proved by  Merino  L{\'o}pez \cite{Mer97} that, for any undirected graph $G$ and any choice of the sink vertex $s$,
\begin{equation*} \label{equation: Merino's theorem for undirected graphs}
 c_0 + c_1y+\ldots+ c_ny^n= y^{|E(G)|} \mathcal T (G;1,y) \qquad \text{(Merino's Theorem)}, 
\end{equation*}
where $\mathcal T(G;x,y)$ is the Tutte polynomial of the graph and $c_i$ is the number of recurrent configurations with $i-\deg(s)$ chips.
(We remark that the extra factor $y^{|E(G)|}$ does not appear in the right side of \cite[Theorem~3.6]{Mer97} as their left side differs from ours by the same  factor.)
As the  Tutte polynomial is defined without any involvement of the vertex $s$,
this implies that 
 $c_0 + c_1y+\ldots+ c_ny^n$  does not depend on the choice of $s$.
 
The sink independence of the polynomial $c_0 + c_1y+\ldots+ c_ny^n$
 was then extended to all Eulerian directed graphs by Perrot and Pham~\cite{Perrot-Pham}, 
 and  they  observed that the same statement does not hold for non-Eulerian directed graphs.
However, they conjectured that a variant of this polynomial has the sink independence property for all directed graphs.   

   Perrot and Pham defined an equivalence relation on the recurrent configurations~(Definition~\ref{definition: sink equivalence relation}),
   and  they defined the total number of chips  of an equivalence class to be the maximum of the  total number of chips of  configurations contained in the class.
   The conjecture of Perrot and Pham is that the sink independence property is true for the polynomial,
   \begin{equation*}
\B(G,s;y):=  c_0' + c_1'y+\ldots+ c_n'y^n,
\end{equation*}
where $c_i$  is the number of equivalence classes with  $i-\outdeg(s)$ chips.
We prove their conjecture by expressing  $\B(G,s;y)$  as an invariant of the sinkless sandpile model (which, as its name implies, does not involve any choice of sink vertex).  
\begin{theorem}[Weak version of Merino's theorem {\cite[Conjecture~1]{Perrot-Pham}}]\label{c. Perrot-Pham} 
  Let $G$ be a strongly connected digraph.
  Then the polynomial $\B(G,s;y)$ is independent of the choice of the vertex $s$. 
\end{theorem}
We call $\B(G,s;y)$  the   \emph{Biggs-Merino polynomial} to honor the contribution of Biggs and Merino  L{\'o}pez to this subject. 

One consequence of Merino's Theorem  is that, for any undirected graph and for any $i$,
the number of recurrent configurations with  
 $i-\deg(s)$ chips is equal to 
the number of spanning trees with external activity $i$.
A bijective proof of this statement was given by Cori and Le Borgne~\cite{CL03},
and we generalize  the bijection of Cori and Le Borgne to all Eulerian directed graphs.

Let $G$ be a strongly connected digraph, and let $s$ be a vertex of $G$.
A \emph{reverse $G$-parking function}~\cite{PS04} with respect to $s$ is a function $f:V(G)\setminus \{s\} \to \N_0$ such that, for any non-empty subset $A \subseteq V(G) \setminus \{s\}$,
there exists $v \in A$ for which $f(v)$ is strictly smaller than the number of edges from $V(G) \setminus A$ to $v$.
An  \emph{arborescence}  of $G$  {rooted} at $s$ is a subgraph of $G$  that contains $|V(G)|-1$ edges and such that for any vertex $v$ of $G$ there exists a unique directed path from $s$ to $v$ in the subgraph. 

We show that Cori-Le Borgne bijection generalizes to a bijection between 
reverse $G$-parking functions and arborescences of $G$ for all directed graphs.
For the full description of the bijection, see  Algorithm~\ref{algorithm: parking functions to arborescences}. 
\begin{theorem}\label{t. generalized Cori-Le Borgne}
Let $G$ be a strongly connected digraph, and let $s \in V(G)$.
Then Cori-Le Borgne bijection generalizes to  a bijection that sends   reverse $G$-parking functions with respect to $s$ to arborescences of $G$ rooted at $s$.
Furthermore, the external activity of the output arborescence is the level of the input reverse $G$-parking function.
\end{theorem}
The  external activity of an arborescence and
the level of a $G$-parking function   are defined in Definition~\ref{definition: external activity} and Definition~\ref{definition: level of a function}, respectively.

As a consequence of Theorem~\ref{t. generalized Cori-Le Borgne} and  
the duality between reverse $G$-parking functions and recurrent configurations for Eulerian directed graphs~\cite[Theorem~4.4]{HLM08},
we get the extension of Merino's Theorem for Eulerian directed graphs.

\begin{theorem}[Merino's Theorem for Eulerian directed graphs] \label{t. Merino's theorem}
Let $G$ be a connected Eulerian digraph.
Then for any $s \in V(G)$,
\[c_0 + c_1y+\ldots+ c_ny^n= t_0+t_1y+\ldots + t_n y^n,      \]
where $c_i$ is the number of recurrent configurations with
 $i-\outdeg(s)$ chips  and $t_i$ is the number of arborescences of $G$ rooted at $s$ with external activity $i$.
 \end{theorem}
The right side of Theorem~\ref{t. main theorem} is known in the literature as the \emph{greedoid polynomial}~\cite{BKL85}, and it can be considered as a single variable generalization of the Tutte polynomial for directed graphs~\cite{GM89,GT90,GM91}.

This paper is arranged as follows:  In Section \ref{s. preliminaries}  we give a review on  the sinkless sandpile model and the sandpile model with a sink.
In Section \ref{s. connections} we prove Theorem~\ref{c. Perrot-Pham}.
In Section~\ref{s. examples} we prove a recurrence relation for the Biggs-Merino polynomial.
In Section \ref{s. greedoid polynomial} we prove  Theorem~\ref{t. generalized Cori-Le Borgne} and Theorem~\ref{t. Merino's theorem}.
 Finally in Section \ref{s. conjecture} we include a list of questions for future 
    research. 
 
\section{Review of abelian sandpile models}\label{s. preliminaries}
In this section we review basic results concerning the sinkless sandpile model and
the sandpile model with a sink.
 We refer to \cite{BL92,HLM08} for a more detailed introduction of these two models, to \cite{PPW13} for an algebraic treatment of this model, 
 and to \cite{BL13,BL142, BL14} for a generalization of these models. 

We use   $G:=(V(G),E(G))$ to denote a \emph{directed graph} (\emph{digraph} for short), possibly with loops and multiple edges.
We use $V$ and $E$ as a shorthand for $V(G)$ and $E(G)$ when the digraph $G$ is evident from the context.
 Each edge $e \in E$ is directed from its source vertex  to its target vertex. 
The  \emph{outdegree} of a vertex $v \in V$, denoted by $\outdeg(v)$, is the number of  edges with $v$ as source vertex,  while the  \emph{indegree} of a vertex $v \in V$, denoted by $\indeg(v)$, is the number of  edges with $v$ as target vertex. 

In this paper, we identify   an  {undirected graph} $G$ with the directed graph obtained by replacing each undirected edge $e:=\{i,j\}$ of $G$ with two directed edges $(i,j)$ and $(j,i)$.
A digraph  obtained in this way is called \emph{bidirected}.


A digraph $G$ is  \emph{Eulerian} if $\outdeg(v)=\indeg(v)$ for all $v \in V$.
In particular, all bidirected graphs are Eulerian.
A digraph  is \emph{strongly connected} if for any two vertices $v,w \in V$ there exists a directed path from $v$ to $w$.
Note  that a connected Eulerian digraph is always strongly connected. 
Throughout this paper, we  always assume that our digraph $G$ is strongly connected.

The  \emph{Laplacian matrix} $\LL$ of a digraph $G$ is the  square matrix $(\LL_{i,j})_{V\times V}$  given by:
\begin{align*}
\LL_{i,j}:=\begin{cases} \outdeg(i)- \# \text{ of loops in vertex } i &\text{ if  }i=j;\\
- \text{ the number of  of edges  from vertex }j \text{ to vertex } i  &\text{ if  }i\neq j. 
\end{cases}
\end{align*}  
Note that the Laplacian matrix is a symmetric matrix if and only if $G$ is bidirected.


\begin{definition}[Primitive period vector]
A vector $\rr \in \R^V$ is a \emph{period vector} of $G$  if it is non-negative, integral, and  $\rr \in \ker(\Delta)$.
A period vector is \emph{primitive} if its entries has no non-trivial common divisor.
\end{definition}
If $G$ is a strongly connected digraph, then 
  a primitive period vector exists, is unique, and is strictly positive~\cite[Proposition~4.1(i)]{BL92}.
Throughout this paper, we will use $\rr$ to denote the primitive period vector of $G$.

 Note that any nonzero period vector of $G$ is a positive multiple of the primitive period vector~\cite[Proposition~4.1(iii)]{BL92}.
Also note  that the primitive period vector  is equal to $(1,\ldots, 1)$ if and only if $G$ is an Eulerian digraph~\cite[Proposition~4.1(ii)]{BL92}.

A \emph{reverse arborescence} of $G$ rooted at $v \in V$ 
 is a subgraph $G$ that contains $|V|-1$ edges and such that for any $w\in V$ there exists a unique directed path from $w$ to $v$ in the subgraph. 
Let $t_v$ denote the number of reverse arborescences of $G$ rooted at $v$.
Note that $t_v>0$ for all $v \in V$ if $G$ is strongly connected.

\begin{definition}[Period constant]\label{definition: period constant}
 The \emph{period constant} $\alpha$ of a strongly connected digraph $G$ is 
\[ \alpha:= \underset{v \in V}{\gcd} \, \{ t_v \}. \qedhere \]
\end{definition}

If $G$ is strongly connected,
then 
by the markov chain tree theorem~\cite{AT89}
the primitive period vector $\rr$ is given by
\[ \rr(v)=\frac{t_v}{\alpha} \quad (v \in V).  \]

\subsection{Sinkless abelian sandpile model}\label{ss. abelian sandpile model}
The  \emph{sinkless (abelian) sandpile model} on a strongly connected digraph $G$, denoted by $\Sand(G)$, 
starts with a number of chips at each vertex of $G$. 
A  \emph{sinkless (chip) configuration} $\cc$  is a vector in $\N_0^{V}$, with $\cc(v)$ representing the number of chips in the  $v \in V$.
A \emph{sinkless firing move} on $\cc$ consists of removing $\outdeg(v)$ chips from a vertex $v$ and sending each of those chips along each edge leaving $v$ to a neighbor of $v$.
Denote by $\satu_v$ the vector in $\N_0^V$ given by $\satu_v(v):=1$ and $\satu_v(w):=0$ for all $w \in V \setminus \{v\}$. 
Note that a sinkless firing move changes  a sinkless configuration $\cc$ to  $\cc-\Delta \satu_v$.

Note that the result of a sinkless firing move is not necessarily a sinkless configuration,
as the $v$-th entry of $\cc-\Delta \satu_v$ 
is negative if $\cc(v)\leq \outdeg(v)$.
A sinkless firing move is \emph{legal} if the fired vertex has at least as many chips as its outgoing edges, or equivalently if the result of the firing move is another sinkless  configuration.
A sinkless configuration $\cc$ is \emph{stable} if $\cc(v)<\outdeg(v)$ for all $v \in V$, or equivalently if there are no legal sinkless firing moves for $\cc$.



Each finite (possibly empty) sequence of sinkless firing moves is associated with 
an  \emph{odometer} $\qq \in \N_0^{V}$, where $\qq(v)$ is equal to the number of times the vertex $v$ is being fired in the sequence.
Note that applying a finite sequence of firing moves with odometer $\qq$ to a sinkless configuration $\cc$ gives us the sinkless configuration $\cc-\Delta\qq$.


%

For any two sinkless configurations $\cc, \dd$ of $G$, 
we write $\cc \dar \dd$  if there exists  a finite (possibly empty) sequence of legal sinkless firing moves that sends  $\cc$ to  $\dd$. 
If the odometer $\qq$ of the sequence  is known,
we will write $\cc \xlongrightarrow[\qq]{} \dd$ instead.

It follows from the definition that $\dar$ is a transitive relation.

\begin{definition}[Recurrent sinkless configurations]
\label{definition: recurrent sinkless}
Let $G$ be a strongly connected digraph.
A sinkless configuration $\cc$ of $G$  is \emph{(sinkless) recurrent} if it satisfies these two conditions:
\begin{itemize}
\item The configuration $\cc$ is not stable; and
\item If $\dd$ is a sinkless configuration that satisfies $\cc \dar \dd$,
then $\dd \dar \cc$. \qedhere
\end{itemize}
\end{definition}
 We  use $\Rec(G)$ to denote the set of all recurrent sinkless configurations of $G$.  



 Recall that  $\rr$ denote  the primitive period vector of $G$.
\begin{lemma}[{\cite[Lemma~1.3, Lemma~4.3]{BL92}}]\label{lemma: BL92}
Let $G$ be a strongly connected digraph.
\begin{enumerate}
\item \label{item: BL92 1} If $\qq_1, \qq_2 \in \N_0^V$ satisfy $\qq_1 \leq \qq_2$ and $\cc,\dd_1,\dd_2$ are sinkless configurations that satisfy $\cc \xlongrightarrow[\qq_1]{} \dd_1$ and  $\cc \xlongrightarrow[\qq_2]{} \dd_2$,
then $\dd_1 \xlongrightarrow[\qq_2-\qq_1]{} \dd_2$.
\item \label{item: BL92 2} If $\cc$ is a sinkless configuration  that satisfies $\cc \xlongrightarrow[k\rr]{} \cc$ for some positive $k$, then $\cc \xlongrightarrow[\rr]{} \cc$. \qed
\end{enumerate}
\end{lemma}



  In the next lemma we  present a  test called the  \emph{sinkless burning test} that checks  whether a given sinkless configuration is recurrent.  
 \begin{proposition}[Sinkless burning test]
 \label{l. burning test for sandg}
 Let $G$ be a strongly connected digraph.
 A sinkless configuration $\cc$ is  recurrent if and only if there exists a finite sequence of legal firing moves from $\cc$ back to $\cc$ such that each vertex $v \in V$ is fired exactly $\rr(v)$ times.
 \end{proposition}
\begin{proof}
Proof for the $\Rightarrow$ direction:  Let $\cc$ be an arbitrary  recurrent sinkless configuration.
Since $\cc$ is not stable by definition of recurrence, 
 there exists a sinkless configuration $\dd$ and a vertex $v \in V$ 
such that $\cc \xlongrightarrow[\satu_v]{} \dd$.
Since $\cc$ is recurrent, there exists $\qq \in \N_0^V$ such that 
$\dd \xlongrightarrow[\qq]{} \cc$.
Write $\qq':=\satu_v+\qq$.
It follows that $\cc \xlongrightarrow[\qq]{} \cc$,
and in particular we have
 $\qq$ is a nonzero period vector of $G$.
 This implies that $\qq$ is a positive multiple of $\rr$,
 and  Lemma~\ref{lemma: BL92}\eqref{item: BL92 2} then implies that $\cc \xlongrightarrow[\rr]{}\cc$, as desired.

Proof for the $\Leftarrow$ direction: 
%
Since $\cc \xlongrightarrow[\rr]{} \cc$ by assumption and $\rr$ is a strictly positive vector, we conclude that $\cc$ is not a stable sinkless configuration.

Let $\dd$ be an arbitrary sinkless configuration that satisfies $\cc \dar \dd$.
It suffices to show that $\dd \dar \cc$.
Let $\qq$ be the odometer of this sequence of firing moves that sends $\cc$ to $\dd$.
Since $\rr$ is a strictly positive vector, there exists a positive $k$ 
such that $\qq \leq k \rr$.

Since $\cc \xlongrightarrow[\rr]{} \cc$ by assumption, it follows that 
$\cc \xlongrightarrow[k\rr]{} \cc$.
On the other hand, we also have $\cc \xlongrightarrow[\qq]{} \dd$ by assumption.
Lemma~\ref{lemma: BL92}\eqref{item: BL92 1} then implies that $\dd \xlongrightarrow[k\rr -\qq]{} \cc$, as desired.
\end{proof}


The next lemma gives two sufficient conditions for a sinkless configuration to be recurrent.

\begin{lemma}\label{p. recurrence can be inherited}
Let $G$ be a strongly connected digraph.
\begin{enumerate}
\item \label{item: recurrence inheritance 1} 
If $\cc$ is a recurrent sinkless configuration and $\dd$ is a sinkless configuration that satisfies $\cc \dar \dd$,
then $\dd$ is a recurrent sinkless configuration.

\item \label{item: recurrence inheritance 2} If $\cc$ is a recurrent sinkless configuration, then for any $k \in \N_0$ and $v \in V$ the sinkless configuration $\cc+k\satu_v$  is also recurrent. 
\end{enumerate}
\end{lemma}

\begin{proof}
\begin{enumerate}
\item We first show that $\dd$ is not a stable sinkless configuration.
Since $\cc \dar \dd$ by assumption and $\cc$ is recurrent, we conclude that $\dd \dar \cc$.
If the odometer of the sequence of legal firing moves from $\dd$ to $\cc$ is nonzero, then $\dd$ is not  stable by definition.
If the odometer is the zero vector, then $\dd=\cc$ is recurrent and hence is not stable.

We now show that if $\dd'$ is a sinkless configuration that satisfies $\dd \dar \dd'$, then $\dd' \dar \dd$.
Since $\cc \dar \dd$ and $\dd \dar \dd'$, the transitivity of $\dar$ implies that $\cc \dar \dd'$.
Since $\cc$ is recurrent, we then have $\dd' \dar \cc$.
Since we also have $\cc \dar \dd$,
the transitivity of $\dar$ then implies that $\dd' \dar \dd$.  
The proof is complete.

\item Since $\cc$ is recurrent, 
we have $\cc \xlongrightarrow[\rr]{} \cc$ by Proposition~\ref{l. burning test for sandg}.
It then follows from the definition of legal firing moves that
$\cc+k\satu_v \xlongrightarrow[\rr]{} \cc+k\satu_v$.
Proposition~\ref{l. burning test for sandg} then implies that $\cc+k\satu_v$ is recurrent, as desired. \qedhere
\end{enumerate}
\end{proof}

\subsection{Abelian sandpile model with a sink}\label{ss. sink}
Let $s \in V$ be a fixed vertex which we refer to as the \emph{sink}.
 The  \emph{(abelian) sandpile model with  a sink} at $s$, denoted by $\Sand(G,s)$, is a variant of the sinkless sandpile model for which 
 the sink vertex $s$ never fires and all chips sent to $s$ are removed from the game.
%

A  \emph{sink (chip) configuration} $\cs$ is  a vector $\N_0^{V}$
such that $\cs(s)=0$.
 A  \emph{sink firing move} 
consists of 
reducing the number of chips of $\cs$ at a vertex $v \in V \setminus \{s\}$ by $\outdeg_G(v)$, and then sending one chip along each outgoing edge of $v$ to its neighbouring vertex that is not $s$.
A sink firing move  is  \emph{legal}  if the fired vertex $v$ has at least as many chips as its outdegree before the firing. 
It is convenient for us  to be able to  fire the sink vertex $s$ as a legal sink firing move, so we adopt the convention that firing $s$ 
is a legal sink firing move that sends a sink configuration 
 $\cs$ back to $\cs$.
 (Note that a legal sink firing move sends a sink configuration to another sink configuration.)

The \emph{odometer} of a sequence of sink firing moves is the vector $\qq \in \N_0^V$ that records the number of times a vertex is fired in the sequence.
Let $\Delta_s$ denote the $V \times V$ matrix obtained by changing the row of the Laplacian matrix $\Delta$ that corresponds to $s$ with the zero vector.
Note that applying a finite sequence of sink firing moves with odometer $\qq$ to a sink configuration $\cs$ gives us the sink configuration $\cs-\LL_s\qq$, provided that $\qq(s)=0$.

For two sink configurations $\cs$ and $\ds$,
we write  $\cs \xlongrightarrow{s} \ds$ if there exists a sequence of legal sink firing moves from $\cs$ to $\ds$.
A sink configuration $\cs$ is  \emph{stable} if   $\cs(v)< \outdeg(v)$ for all $v \in V$.

\begin{definition}
[Stabilization]
For any  sink configuration $\cs$,
the \emph{stabilization} $\cs^\circ$ of $\cs$ is a sink configuration such that  $\cs \xlongrightarrow{s} \cs^\circ$ and  $\cs$ is a stable sink configuration. \qedhere
\end{definition}
For a strongly connected digraph $G$,
any sink configuration $\cs$ has a unique stabilization~\cite[Lemma~2.4]{HLM08}.

\begin{lemma}\textnormal{(\cite[Lemma~2.2]{HLM08}).}\label{l. e and f related get same stabilizer}
Let $G$ be a strongly connected digraph, let $s\in V$, and let $\cs$ and $\ds$ be sink configurations of $G$.
If $\cs \xlongrightarrow{s} \ds$, then $\cs^\circ=\ds^\circ$. \qed
\end{lemma}
%

\begin{definition}[Recurrent sink configurations]
\label{definition: recurrent sink}
Let $G$ be a strongly connected digraph, and let $s\in V$.
A sink configuration $\cs$ of $G$ is \emph{(sink) recurrent} if for any sink configurations $\cs_1$ there exists another sink configuration $\cs_2$ such that $(\cs_1+\cs_2)^\circ=\cs$.
\end{definition}
Note that a recurrent sink configuration is always a stable configuration.
We  use $\Rec(G,s)$ to denote the set of s-recurrent configurations of $\Sand(G,s)$. 
When there is a possible ambiguity between the two notion of recurrence,
 \emph{sinkless recurrence} will refer to Definition~\ref{definition: recurrent sinkless}, 
and  \emph{sink recurrence} will refer to Definition~\ref{definition: recurrent sink}.


%
\begin{lemma}[Abelian property \textnormal{\cite[Corollary~2.6]{HLM08}}]\label{l. abelian property}
Let $G$ be a strongly connected digraph, let $s \in V$, and let 
  $\cs_1, \cs_2,\cs_3$ be sink  configurations.
   Then:
    \[
\pushQED{\qed}     
    ((\cs_1+\cs_2)^\circ+\cs_3)^\circ=((\cs_1+\cs_3)^\circ+\cs_2)^\circ =(\cs_1+\cs_2+\cs_3)^\circ.  
    \popQED
    \] 
\end{lemma}


In the next proposition we  present a burning test to check for sink recurrence.
It is  first discovered by Dhar~\cite{Dhar90} for undirected graphs and then  by Speer~\cite{Speer93} and Asadi and Backman~\cite{AB11} for directed graphs. 

The \emph{sink Laplacian vector} $\us \in \N_0^V$ is 
\[
\us(v):=\begin{cases} \text{number of edges from $s$ to $v$ } &\text{ if  }v\neq s;\\
0 &\text{ if  }v=s.  
\end{cases}   
\]
Recall that $\rr(s)$ is the entry of the primitive period vector $\rr$ that corresponds to $s$.


\begin{proposition}[Sink burning test  \textnormal{\cite[Theorem~3]{Speer93}},  \textnormal{\cite[Theorem~3.11]{AB11})}]\label{l. burning test thief} 
Let $G$ be a strongly connected digraph, let $s \in V$, and let $\cs$ be a sink configuration of $G$.
Then $\cs$ is  a sink recurrent configuration if and only if 
 $(\cs+ \rr(s) \us)^\circ=\cs$.  \qed 
\end{proposition}


The next lemma gives a sufficient condition for a sink configuration to be recurrent.
\begin{lemma}\textnormal{(\cite[Lemma~2.17]{HLM08}).}\label{p. s-recurrence can be inherited}
Let $G$ be a strongly connected digraph, let $s \in V$, and let $\cs$ be a sink recurrent configuration.
If $\ds$ is a sink configuration such that there exists a sink configuration 
$\cs'$ satisfying $\ds=(\cs+\cs')^\circ$, then
 $\ds$ is a sink recurrent configuration. \qed
\end{lemma}

Let $Z_s \subseteq \Z^V$ denote the set
\[  Z_s:=\{ \zz \in \Z^V \mid \zz(s)=0  \}.\]
Note that $\Rec(G,s)$  is a subset of $Z_s$.
Also note that $\Delta_s Z_s \subseteq Z_s$.

\begin{lemma}[{\cite[Corollary~2.16, Corollary~2.18]{HLM08}}]\label{lemma: sandpile group bijection}
Let $G$ be a strongly connected digraph, and let $s \in V$.
Then
\begin{enumerate}
\item \label{item: sandpile group 1}  The inclusion map $\Rec(G,s) \to Z_s/\Delta_s Z_s$ is a bijection. 
\item \label{item: sandpile group 2} The cardinality of $\Rec(G,s)$ is equal to the number of reverse arborescences of $G$ rooted at $s$.  \qed
\end{enumerate}
\end{lemma}

\subsection{A Connection between the sinkless sandpile model and the sandpile model with sink}\label{ss. comparion two sandpile models}

Let $\cc$ be a sinkless configuration of $G$.
In order to reduce the number of notations, 
we denote by $\cs$ the sink configuration of $G$
with $\cs(v):=\cc(v)$ if $v \neq s$ and $\cs(v):=0$  if $v=s$.


Let $v_1,\ldots, v_k$ be vertices of $G$.
For two sinkless configurations $\cc$ and $\dd$,
we write $\cc \xlongrightarrow[v_1\cdots v_k]{}\dd$
the sequence of sinkless firing moves that fires $v_1,\ldots, v_k$ (in that order) is legal and sends $\cc$ to $\dd$.
For two sink configurations $\cs$ and $\ds$,
we write $\cs \xlongrightarrow[v_1\cdots v_k]{s}\ds$
the sequence of sink firing moves that fires $v_1,\ldots, v_k$ (in that order) is legal and sends $\cs$ to $\ds$.

In the next lemma we highlight a connection between sinkless configurations and sink configurations.


\begin{lemma}\label{p. sequence of firing moves is transferrable}
Let $\cc$ and $\dd$ be sinkless configurations.
\begin{enumerate} 
\item \label{item: tranferrable 1} 
Let $v_1,\ldots,v_k \in V \setminus \{s\}$, 
and let $n$ be the number of chips removed from the game by the sequence of sink firing moves that fires $v_1,\ldots, v_k$.
If 
 $\cs \xlongrightarrow[v_1\cdots v_k]{s}\ds$ and $\dd(s)-\cc(s)=n$,
then   $\cc \xlongrightarrow[v_1\cdots v_k]{}\dd$.
\item \label{item: transferrable 2}
Let $v_1,\ldots, v_k \in V$, and let   $m$ be the number of instances of $s$ in the sequence $v_1,\ldots, v_k$.
 If 
 $\cc \xlongrightarrow[v_1\cdots v_k]{}\dd$,
 then 
 $\cs+m \us \xlongrightarrow[v_1\cdots v_k]{s}\ds$.
\end{enumerate}
\end{lemma}
\begin{proof}
\begin{enumerate} 
\item By induction on $k$, it suffices to prove the claim for when $k=1$.
Since firing $v_1$ is a legal sink firing move on $\cs$ and $v_1\neq s$,
we have firing $v_1$ is also a legal sinkless firing move on $\cc$.
Now note that \begin{align*}
\cc- \Delta \satu_{v_1}=&\cc +n\satu_s -\Delta_s \satu_{v_1}  \quad \text{(since }v_1 \neq s)\\
=&\dd+ (\cc- \dd) +n\satu_s -\Delta_s \satu_{v_1} = \dd + (\cs- \ds)  -\Delta_s \satu_{v_1}=\dd.
\end{align*}
Hence we conclude that $\cc \xlongrightarrow[v_1]{}\dd$, as desired.

%

\item By induction on $k$, it suffices to prove the claim for when $k=1$.

First consider the case when $v_1 =s$.
Note that by definition the legal sink firing move that fires $s$ 
sends $\cs+\us$ back to $\cs+\us$.
On the other hand, we have 
\begin{align*}
\dd= \cc- \Delta \satu_s= \cc+ \us- \outdeg(s)\satu_s.
\end{align*}
This then implies that $\ds=\cs+\us$.
Hence we have  $\cs+m \us \xlongrightarrow[v_1]{s}\ds$.

Now consider the case when $v_1\neq s$.
Since firing $v_1$ is a legal sinkless firing move on $\cc$
and $v_1\neq s$,
we have firing $v_1$ is also a legal sink firing move on $\cs$.
Now note that
\begin{align*}
\dd=\cc-\Delta_s \satu_{v_1}=&\cc +n\satu_s -\Delta_s \satu_{v_1}  \quad \text{(since }v_1 \neq s).
\end{align*}
This implies that $\ds=\cs-\Delta_s \satu_{v_1}$.
Since $m=0$ when $v_1\neq s$,
we conclude that 
 $\cs+m \us \xlongrightarrow[v_1]{s}\ds$.
 The proof is complete.
\qedhere
\end{enumerate}
\end{proof}

\section{Proof of Theorem~\ref{c. Perrot-Pham}}\label{s. connections}

In this section we prove the conjecture of Perrot and Pham~\cite[Conjecture~1]{Perrot-Pham}.


Recall that $\Delta$ is the Laplacian matrix of $G$.
\begin{definition}[Sinkless equivalence relation]
For any recurrent sinkless configurations $\cc$ and $\dd$,
we write $\cc \sim \dd$ if there exists $\zz \in \Z^V$ such that $\cc-\dd=\Delta \zz$.
\end{definition}
Note that $\sim$ defines an equivalence relation on the set of recurrent sinkless configurations.
We call an equivalence class for the relation $\sim$  a \emph{recurrent (sinkless) class}.
For any recurrent sinkless configuration $\cc$,
we denote by $[\cc]$ the recurrent sinkless class  that 
contains $\cc$.

\begin{definition}
[Sinkless level]
For any sinkless configuration $\cc$, the \emph{level} of $\cc$, denoted by $\lvl(\cc)$, is the total number of chips in the configuration $\cc$.
For any recurrent sinkless class $[\cc]$, the \emph{level} of $[\cc]$, denoted by $\lvl([\cc])$, is the level of a configuration  contained in $[\cc]$. 
\end{definition}
It is straightforward to check that two sinkless configurations has the same total number of chips if they are related by $\sim$, and hence  the level of recurrent sinkless classes is well defined.

Recall that $\Delta_s$ is the $V \times V$ matrix obtained by changing the row of the Laplacian matrix of $G$ that corresponds to $s$ with the zero vector.
\begin{definition}[Sink equivalence relation]\label{definition: sink equivalence relation}
For any  recurrent sink configurations $\cs$ and $\ds$,
we write $\cs \es \ds$ if there exists $\zz \in \Z^V$ such that $\cs-\ds=\Delta_s \zz$.
\end{definition}
Note that $\es$ defines an equivalence relation on the set of recurrent sink configurations.
We call an equivalence class for the relation $\es$  a \emph{recurrent (sink) class}. 
For any recurrent sink configuration $\cs$,
we denote by $[\cs]_s$ the recurrent sink class  that 
contains $\cs$.

\begin{definition}
[Sink level]
For any sink configuration $\cs$, the \emph{level} $\lvl(\cs)$ of $\cs$ is the total number of chips in $\cs$.
For any recurrent sink class $[\cc_s]$, the \emph{level} of $[\cs]_s$ is 
\[\lvl([\cs]_s):=\max \{ \lvl(\ds) \mid \ds \in [\cs]_s   \}. \qedhere \] 
\end{definition}

For any nonnegative integer $m$,
we denote by $\Rec_m(G,\sim)$ the set of recurrent sinkless classes with level $m$,
and by $\Rec_m(G,\es)$ the set of recurrent sink classes with level $m$.


 \begin{proposition}\label{p. size Rec(G,es)}
 Let $G$ be a strongly connected digraph, and let $s \in V$.
 Then 
 the cardinality of  $\Rec(G,\es)$ is equal to the period constant $\alpha$ of $G$. 
 \end{proposition}
\begin{proof}
It follows from Lemma~\ref{lemma: sandpile group bijection}\eqref{item: sandpile group 1} and the definition of $\es$ that:
\begin{align*}
|\Rec(G,\es)|=\left| \frac{Z_s}{\Delta_s \Z^{V}} \right|,
\end{align*}
where $Z_s$ and $\Delta_s$ is as defined in Section~\ref{s. preliminaries}.
Now note that 
\begin{align*}
\left| \frac{Z_s}{\Delta_s \Z^{V}} \right|&={\left| \frac{Z_s}{\Delta_s Z_s} \right|} \, \bigg/ \, {\left| \frac{\Delta_s \Z^V }{\Delta_s Z_s} \right|} \quad \text{(by the third isomorphism theorem for groups)}\\
&= {t_s} \, \bigg/\,  { \left| \frac{\Delta_s \Z^V}{\Delta_s Z_s}  \right|} \quad \text{(by Lemma~\ref{lemma: sandpile group bijection}\eqref{item: sandpile group 2})}.
\end{align*}
 By a direct computation, we have
 $|{\Delta_s \Z^V} \, / \, {\Delta_s Z_s} |$ is equal to $\rr(s)$,
 where $\rr$ is the primitive period vector of $G$.
Now recall that $\rr(s)$ is equal to $t_s/\alpha$ by the markov chain tree theorem~\cite{AT89}. 
 Hence we conclude that
 \begin{align*}
|\Rec(G,\es)|={t_s} \, \bigg/\,  { \left| \frac{\Delta_s Z_s}{\Delta_s \Z^{V}}  \right|} =\frac{t_s}{\rr(s)}=\alpha,
\end{align*}
as desired.
\end{proof}
As a corollary of Proposition~\ref{p. size Rec(G,es)}, 
we have the cardinality of $\Rec(G,\es)$ is independent of the choice of $s$.
We will prove   a  stronger  sink independence result  in the next theorem.

%

\begin{definition}[Biggs-Merino polynomial]
Let $G$ be a strongly connected digraph, and let $s \in V$.
The \emph{Biggs-Merino polynomial}
  $\B(G,s;y)$ is 
\[\B(G,s;y):= \sum_{m\geq 0} |\Rec_m(G,\es)|\cdot y^{m+\outdeg(s)}. \qedhere
\]
\end{definition}
Since $\Rec(G,\es)$ is a finite set by Proposition~\ref{p. size Rec(G,es)}, 
we have $\Rec_m(G,\es)$ is an empty set  for sufficiently large $m$.
This then implies that $\B(G,s;y)$ is a polynomial.

We denote by $\Rc(G;y)$  the formal power series 
\[\Rc(G;y):= \sum_{m\geq 0} |\Rec_m(G,\el)|\,  y^{m}.\]

\begin{theorem}\label{t. main theorem}
Let $G$ be a strongly connected digraph, and let $s \in V$.
 We have the following equality of formal power series:
\begin{equation*}
\Rc(G;y)=\frac{\B(G,s;y)}{(1-y)}.
\end{equation*}
\end{theorem}

The following conjecture of Perrot and Pham~\cite{Perrot-Pham} is a direct corollary of Theorem~\ref{t. main theorem}.

  \begin{reptheorem}{c. Perrot-Pham} \textnormal{ (\cite[Conjecture~1]{Perrot-Pham}).}
  Let $G$ be a strongly connected digraph.
  Then $\B(G,s;y)$ is independent of the choice of the vertex $s$. \qed
  \end{reptheorem}


%


The rest of this section is focused on the proof of Theorem~\ref{t. main theorem}.

For any nonnegative $m$, we define the map $\varphi$ by
\begin{align*}
\varphi: \Rec_m(G, \sim) &\to \bigsqcup_{n \leq m-\outdeg(s)}\Rec_{n } (G, \es)\\
[\cc] &\mapsto [\cs^\circ]_s
\end{align*}

The following lemma shows that $\varphi$ is well defined and is injective for all positive $m$.
\begin{lemma}\label{lemma: varphi is well defined and injective}
Let $G$ be a strongly connected digraph, let $s \in V$,
and let $\cc,\dd$ be sinkless recurrent configurations of $G$ with the same level.
Then
\begin{enumerate}
\item  The sink configuration $\cs^\circ$ is  sink recurrent,  and  $\lvl(\cs^\circ) \leq \lvl(\cc)-\outdeg(s)$.
\item $\cc \sim \dd$ if and only if $
\cs^\circ \es \ds^\circ$.
\end{enumerate}
\end{lemma}
\begin{proof}
\begin{enumerate}
\item 
Since $\cc$ is sinkless recurrent, 
we have $\cc \xlongrightarrow[\rr]{} \cc$ 
by  Lemma \ref{l. burning test for sandg} (recall that $\rr$ is the primitive period vector of $G$).
By Lemma \ref{p. sequence of firing moves is transferrable}\eqref{item: transferrable 2},  we then have  $\cs+\rr(s)\us \xlongrightarrow[\rr]{s}\cs$.
By Lemma \ref{l. e and f related get same stabilizer}, this implies that $(\cs+\rr(s)\us)^\circ = \cs^\circ$.
By Lemma  \ref{l. abelian property},  we then conclude that
$(\cs^\circ+\rr(s)\us)^\circ  = (\cs+\rr(s)\us)^\circ = \cs^\circ$.
 Hence   $\cs^\circ$ passes the sink burning test in  Proposition \ref{l. burning test thief},  and we have $\cs^\circ$ is a recurrent sink configuration.

Let $n$ be the number of chips removed during the stabilization of $\cs$, and let
 $\dd:=\cs^\circ+(n+\cc(s))\satu_s$.
By Lemma~\ref{p. sequence of firing moves is transferrable}\eqref{item: tranferrable 1},
we conclude that 
$\cc \xlongrightarrow{} \dd$.
This implies that $\lvl(\cc)=\lvl(\dd)$ as legal firing moves do not change the total number of chips.
By  Lemma~\ref{p. recurrence can be inherited}\eqref{item: recurrence inheritance 1},
this also implies that 
 $\dd$ is a recurrent sinkless configuration.
In particular, we have $\dd$ is not a stable sinkless configuration.

Since $\cs^\circ$ is a recurrent sink configuration (and hence stable),
we have $\cs^\circ(v)<\outdeg(v)$ for all $v \in V$.
This implies that $\dd(v)=\cs^\circ(v) < \outdeg(v)$ for all $v \in V \setminus\{s\}$.
 Since $\dd$ is not a stable sinkless configuration,
 we then conclude that $\dd(s)\geq \outdeg(s)$.
 Now note that
 \[ \lvl(\cc)=\lvl(\dd)=\lvl(\cs^\circ)+\dd(s)\geq \lvl(\cs^\circ)+\outdeg(s),\]
and the proof is complete.

\item Let $\qq_1$  be the odometer of a sequence of sink firing moves that stabilizes $\cs$, and let $\qq_2$  be the odometer of a sequence of sink firing moves that stabilizes $\ds$.
We have
\begin{equation}\label{equation: varphi}
\cc-\dd = \cs-\ds +(\cc(s)-\dd(s))\satu_s =\cs^\circ-\ds^\circ + \Delta_s(\qq_1-\qq_2)+(\cc(s)-\dd(s))\satu_s.
\end{equation}

If $\cc-\dd=\Delta \zz$ for some $\zz \in \Z^V$, then :
\begin{align*}
\cs^\circ-\ds^\circ =& \Delta\zz- \Delta_s (\qq_1-\qq_2)-(\cc(s)-\dd(s))\satu_s \quad \text{(by equation~\eqref{equation: varphi})}\\
=&\Delta_s( \zz-\qq_1+\qq_2) +t\satu_s,
\end{align*}
for some $t \in \Z$.
Since $\cs^\circ(s)=\ds^\circ(s)=0$ and $\satu_s^\top \Delta_s=(0,\ldots,0)$, 
\[0 =\satu_s^\top(\cs^\circ -\ds^\circ)= \satu_s^\top(\Delta_s( \zz-\qq_1+\qq_2) +t\satu_s)=t. \]
Hence we have $\cs^\circ-\ds^\circ=\Delta_s( \zz-\qq_1+\qq_2)$, which implies
 that $\cs^\circ \es \ds^\circ$.

If $\cs^\circ-\ds^\circ=\Delta_s \zz$ for some $\zz \in \Z^V$, then
\begin{align*}
\cc-\dd =& \Delta_s \zz+ \Delta_s (\qq_1-\qq_2)+(\cc(s)-\dd(s))\satu_s \quad \text{(by equation~\eqref{equation: varphi})}\\
=&\Delta( \zz+\qq_1-\qq_2) +t\satu_s,
\end{align*}
for some $t \in \Z$.
Since  $\lvl(\cc)=\lvl(\dd)$ by assumption and $(1,\ldots,1)^\top \Delta=(0,\ldots,0)$, 
\begin{align*}
0=&(1,\ldots, 1)^\top(\cc-\dd) = (1,\ldots,1)^\top (\Delta( \zz+\qq_1-\qq_2) +t\satu_s) 
=t. 
\end{align*}
Hence we have $\cc-\dd =\Delta( \zz+\qq_1-\qq_2)$, which implies that $\cc \sim \dd$.
The proof is complete. \qedhere
\end{enumerate}
\end{proof}

We now proceed by showing that the map $\varphi$ is surjective,
and we need the following technical lemma.

\begin{lemma}\label{lemma: varphi is surjective}
Let $G$ be a strongly connected digraph, let $s \in V$,  and let $\cs$ be a recurrent sink configuration such that 
$\lvl(\cs)\geq \lvl(\ds)$ for all $\ds \in [\cs]_s$.
Let $\cc:=\cs+\outdeg(s) \satu_s$.
Then $\cc$ is a recurrent sinkless configuration of $G$.
\end{lemma}
\begin{proof}

Let $v_1,\ldots, v_k$ be a sequence of legal sinkless firing moves on $\cc$ 
with odometer $\qq$  less than the primitive period vector $\rr$.
Without loss of generality,
assume that $v_1,\ldots, v_k$ is of maximum length.
Note that $k\geq 1$ as firing $s$ is a legal sinkless firing move on $\cc$,
and in particular $\qq$ is a nonzero vector.


Write $\cc':=\cc-\Delta \qq$.
We claim that $\cc'(v)<\outdeg(v)$ for all $v \in V \setminus \{s\}$.
Suppose to the contrary that $\cc'(v)\geq \outdeg(v)$ for some $v \in V \setminus \{s\}$.
By the maximality of the odometer $\qq$,
it follows that $\qq(v) =\rr(v)$.
Now note that
\begin{align*}
\cc'(v)=&\cc(v)+\Delta \qq(v)= \cc(v)+ \outdeg(v) \qq(v) -\sum_{w \in V \setminus \{v\}} \Delta_{v,w} \qq(w)\\
\leq &\cc(v)+ \outdeg(v) \qq(v) -\sum_{w \in V \setminus \{v\}} \Delta_{v,w} \rr(w)\\
= &\cc(v)+ \outdeg(v) \rr(v) -\sum_{w \in V \setminus \{v\}} \Delta_{v,w} \rr(w)=\cc(v)=\cs(v).
\end{align*}
Since $\cs$ is a recurrent sink configuration (and hence stable) and $v \in V\setminus\{s\}$,
we have $\cs(v)<\outdeg(v)$.
This means that $\cs'(v)=\cs(v)<\outdeg(v)$, and we get a contradiction.
This proves the claim.
 
 Since $\cc \xlongrightarrow[\qq]{}\cc'$, 
 we have $\cs+\qq(s) \us  \xlongrightarrow[\qq]{s} \cs'$ by Lemma~\ref{p. sequence of firing moves is transferrable}\eqref{item: transferrable 2}.
 Since $\cs'$ is a stable sink configuration,
 this implies that  $(\cs+\qq(s) \us)^\circ=\cs'$.
 Since $\cs$ is a recurrent sink configuration,
 we have $\cs'$ is a recurrent sink configuration
  by Lemma~\ref{p. s-recurrence can be inherited}.
  It then follows that $\cs'$ is contained in $[\cs]_s$, and hence
  we have $\lvl(\cs)\geq \lvl(\cs')$ by assumption. 
  On the other hand, we have $\lvl(\cc)=\lvl(\cc')$ since $\cc \xlongrightarrow[\qq]{} \cc'$.
  Now note that:
  \begin{align}\label{equation: surjective 1}
  \begin{split}
  0=&\lvl(\cc)-\lvl(\cc')=\lvl(\cs)-\lvl(\cs')+\cc(s)-\cc'(s)\\
  \geq& \cc(s)-\cc'(s)=\outdeg(s) -\cc'(s).
  \end{split}
  \end{align}
Hence we have $\cc'(s)\geq  \outdeg(s)$.
By the maximality of the odometer $\qq$,
this then implies that $\qq(s)=\rr(s)$.

Now note that 
\[\cs'=(\cs+\qq(s) \us)^\circ=(\cs+\rr(s) \us)^\circ=\cs,\]
where the last equality is due to Proposition~\ref{l. burning test thief}.
This implies that we have equality in equation~\eqref{equation: surjective 1},
which then implies that $\cc'(s)=\cc(s)$.
Hence we conclude that $\cc'=\cs'+\cc'(s)=\cs+\cc(s)=\cc$.

Now note that $\qq \in \ker(\Delta)$ since $\Delta \qq= \cc-\cc'=(0,\ldots,0)^\top$.
Since $\ker(\Delta)$ has dimension 1 (as $G$ is strongly connected) and $\qq$ is nonnegative, we conclude that $\qq=k\rr$ for some nonnegative $k$.
Since we have previously shown that $\qq$ is a nonzero vector,
we have that $k$ is positive.
By Lemma~\ref{lemma: BL92}\eqref{item: BL92 2}, we conclude that $\cc \xlongrightarrow[\rr]{} \cc$.
It then follows from Proposition~\ref{l. burning test for sandg} that 
$\cc$ is a recurrent sinkless configuration.
\end{proof}

\begin{lemma}\label{lemma: varphi is a bijection}
Let $G$ be a strongly connected digraph, let $s \in V$, 
and let $m$ be a nonnegative integer.
The map $\varphi: \Rec_m(G, \sim) \to \bigsqcup_{n \leq m-\outdeg(s)}\Rec_{n } (G, \es)$ is a bijection.
\end{lemma}

\begin{proof}
All other properties except for the surjectivity of $\varphi$ have been checked in Lemma~\ref{lemma: varphi is well defined and injective}.

Let $[\cs]_s$ be a recurrent sink class of $G$ with level at most $m-\outdeg(s)$,
and without loss of generality let $\cs$ be a recurrent sink configuration in $[\cs]_s$
such that $\lvl(\cs)\geq \lvl(\ds)$ for all $\ds \in [\cs]_s$.
Let $n:=\lvl(\cs)=\lvl([\cs]_s)$,
and let  $\cc:=\cs+(m-n)\satu_s$.
Note that $\lvl(\cc)=\lvl(\cs)+m-n=m$, and therefore the surjectivity of $\varphi$ follows if we can show 
 that $\cc$ is a recurrent sinkless configuration.

Let $\cc':=\cs+\outdeg(s)\satu_s$.
By Lemma~\ref{lemma: varphi is surjective},
we have $\cc'$ is a recurrent sinkless configuration.
Now note that 
\[\cc(s)=m-n=m-\lvl([\cs]_s)\geq \outdeg(s),\]
and hence $\cc'=\cc+k\satu_s$ for some nonnegative $k$.
It then follows from Lemma~\ref{p. recurrence can be inherited}\eqref{item: recurrence inheritance 2} that $\cc$ is a recurrent sinkless configuration.
The proof is complete.
\end{proof}

\begin{proof}[Proof of Theorem~\ref{t. main theorem}]
By Lemma \ref{lemma: varphi is a bijection},
we have for any nonnegative $m$ 
\begin{equation}\label{equation: main theorem 1}
|\Rec_{m}(G,\el)|= \sum_{n=0}^{m-\outdeg(s)} |\Rec_{n}(G,\es)|.
\end{equation}
Now note that 
\begin{align*}
\frac{\B(G,s;y)}{(1-y)}&=\frac{\sum_{n\geq 0} |\Rec_n(G,\es)| \cdot y^{n+\outdeg(s)}}{(1-y)} \\
&=\left( \sum_{n\geq 0} |\Rec_n(G,\es)| \cdot y^{n+\outdeg(s)}\right) \cdot \left( \sum_{k\geq 0} y^k\right)\\
&=\sum_{m\geq 0} \left(\sum_{n=0}^{m-\outdeg(s)} |\Rec_{n}(G,\es)| \right)y^m=\sum_{m \geq 0} |\Rec_m(G,\el)| y^{n}\\
&=\Rc(G;y). \qedhere
\end{align*}
\end{proof}

\section{A recurrence relation for the Biggs-Merino polynomial}\label{s. examples}
In this section we present a recurrence relation for the Biggs-Merino polynomial, 
and we apply it to
  compute  the  Biggs-Merino polynomial
 for  a(n infinite) family of non-Eulerian digraphs.

\begin{proposition}\label{l. computation}
Let $G$ be a strongly connected digraph, let $s \in V$,  and  let $k$ be a positive natural number.
 Let $G^k$ be   the digraph  obtained from $G$ by replacing 
 each edge in $G$ with $k$ copies  of the same edge.
Then for any $s \in V$,
\[\B(G^k,s;y)= \B(G,s;y^k) \, \left( \frac{1-y^k}{1-y} \right)^{|V|-1}  .\]
\end{proposition}

\begin{proof}
For any sinkless configuration $\cc$ in $G^d$, let $\pi(\cc)$ be the sinkless configuration in $\Sand(G)$ given by
\[  \pi(\cc)(v):=\lfloor {\cc(v)}/{k}\rfloor \quad (v \in V). \] 

%
%

Let $\cc$ and $\cc'$ be configurations in $\Sand(G^k)$.
It is straightforward to check that
 $\cc'$ is accessible from $\cc$ by a sequence of (legal) sinkless firing moves in $G^d$ if and only if 
$\pi(\cc')$ is accessible from $\pi(\cc)$ by the same sequence of (legal) firing moves in $G$ 
and for all $v \in V$ we have $\cc(v)\equiv \cc'(v)$  (mod $k$).

Let $[\dd]$ be a sinkless recurrent class in $G$.
For any $\hh \in \{0,\ldots, k-1\}^V$, let $\dd_\hh$ be the configuration in $\Sand(G^d)$ given by 
\[ \dd_{\hh}(v):= k \, \dd(v)+\hh(v) \quad (v \in V).  \]
It then follows from the conclusion in the previous paragraph that:
\[ \pi^{-1}([\dd])= \bigsqcup_{ \hh \in \{0,\ldots, k-1 \}^{V}  } [\dd_\hh]. \] 

Hence we have
\begin{align*}
\Rc(G^k;y)=\sum_{[\cc] \in \Rec(G^k,\el)} y^{\lvl([\cc])}&= \sum_{[\dd] \in \Rec(G,\el)} \,  \sum_{\hh \in \{0,\ldots, k-1 \}^{V}} y^{\lvl[\dd_\hh]}\\
&=  \sum_{[\dd] \in \Rec(G,\el)} y^{k \, \lvl([\ee])} \,  (1+y+\ldots y^{k-1})^{|V|}  \\
&=      \Rc(G;y^k) \,  \left( \frac{1-y^k}{1-y} \right)^{|V|}.
\end{align*}
Together with Theorem~\ref{t. main theorem}, this implies that
\begin{align*}
\B(G^k,s;y)=&\Rc(G^k;y)\, (1-y)={\Rc(G;y^k)} \, \frac{(1-y^k)^{|V|}}{(1-y)^{|V|-1}}\\
 =&{\B(G,s;y^k)} \, \left( \frac{1-y^k}{1-y} \right)^{|V|-1}. \qedhere
\end{align*}
 \end{proof}

 Let $n,a,b$ be positive integers.
  We denote by  $G(n;a,b)$  the digraph with vertex set  $\{v_1,v_2, \ldots, v_n\}$, and  with $a$ edges from $v_i$ to $v_{i+1}$ 
  and $b$ edges from $v_{i+1}$ to $v_{i}$ for $1 \leq i \leq n-1$.
Note that $G(n;a,b,)$ is  Eulerian if and only if  $a= b$.


\begin{figure}[ht!]
  \centering
    \includegraphics[scale=1.0]{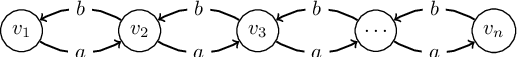}
      \caption{The digraph $G(n;a,b)$ with vertex set $\{v_1,\ldots, v_n\}$, and  with $a$ edges from $v_i$ to $v_{i+1}$ and $b$ edges from $v_{i+1}$ to $v_{i}$ ($i \in \{1,\ldots, n-1\}$).}
\end{figure}

\begin{lemma}\label{lemma: Gnab}
Let $n,a,b$ be positive integers, and let $k:=\gcd(a,b)$.
Then 
\begin{align*}
\B(G(n;a,b),v_1;y)
&= y^{(n-1)(a+b-d)}  \left( \frac{1-y^k}{1-y}\right)^{n-1}. 
\end{align*}
\end{lemma}
\begin{proof}
We start with the case when $\gcd(a,b)=1$.
 Note that, for any $i\in \{1,\ldots, n-1\}$, the number of reverse arborescences of $G(n;a,b)$  rooted at $v_i$ is equal to $b^{n-i}a^{i-1}$.
 Hence the period constant $\alpha$ of $G(n;a,b)$  is equal to
 \[ \alpha= \gcd_{1 \leq i \leq n-1} b^{n-i}a^{i-1} =\gcd(a,b)^{n-1}=1. \]
By Lemma~\ref{p. size Rec(G,es)}, 
this implies that 
  $\Rec(G(n;a,b,),\stackrel{v_1}{\sim})$ contains only one element.


Let $\cs$ be the stable sink configuration of $G(n;a,b)$ with maximum level, i.e.
\begin{align*}
\cs(v):=\begin{cases}
\outdeg(v)-1  & \text{if } v \neq v_1;\\
0  & \text{if } v = v_1.
\end{cases}
\end{align*}

Since $\Rec(G(n;a,b,),\stackrel{v_1}{\sim})$ contains only one element,
we conclude that $[\cs]_{v_1}$ is the unique element in $\Rec(G(n;a,b,),\stackrel{v_1}{\sim})$.
Hence we have:
\begin{align}
\label{equation: Gnab 1}
\begin{split}
 \B(G(n;a,b), v_1;y)=&  y^{\lvl([\cs]_{v_1}) +\outdeg(v_1)}\\
 =&y^{{\lvl(\cs)}+\outdeg(v_1)} \quad \text{(by the maximality of }\cs)\\
=&y^{(n-1)(a+b-1)}. 
\end{split}
\end{align}

We now proceed with the case when $k=\gcd(a,b)$ is arbitrary.
Note that
\begin{align*}
\B(G(n;a,b),v_1;y)&=\B({G(n;a/k,b/k)},v_1;y^k)  \left(\frac{1-y^k}{1-y} \right)^{n-1} \quad \text{(by Proposition~\ref{l. computation})}\\ 
&= y^{(n-1)(a+b-k)}  \left( \frac{1-y^k}{1-y}\right)^{n-1} \quad \text{(by equation~\eqref{equation: Gnab 1})}. \qedhere 
\end{align*}
\end{proof}


By a similar argument as in Lemma~\ref{lemma: Gnab}, for any $k \geq 1$ and any strongly connected digraph $G$ with the period  constant equal to 1,
\[ \B(G^k,s;y)=y^{k(|E(G)|-|V|+1))} \left( \frac{1-y^k}{1-y}\right)^{|V|-1}.\]

\section{Connections to the greedoid  polynomial}\label{s. greedoid polynomial}
In this section we relate the Biggs-Merino polynomial to another invariant of digraphs called the greedoid  polynomial.

\subsection{Greedoid polynomial and reverse G-parking functions}\label{subsection: definition of greedoid polynomial and reverse G-parking function}


%


Let $G$ be a  directed graph.
A  \emph{directed path} $P$ of $G$ of length $k$ is a sequence $e_1\ldots e_{k}$
such that for $i \leq \{1,\ldots, k-1\}$ the target vertex of $e_{i}$ is the source vertex of $e_{i+1}$.

\begin{definition}[Arborescences]
Let $G$ be a strongly connected digraph.
An  \emph{arborescence} $T$ of $G$  {rooted} at $s \in V$ is a subgraph of $G$  that contains $|V|-1$ edges and such that for any $v \in V$ there exists a unique directed path from $s$ to $v$ in the subgraph. 
\end{definition}

Fix a total order $<$ on the directed edges of $G$.
For any two distinct edge-disjoint directed paths $P_1$ and $P_2$,
we write $P_1 < P_2$ if 
the smallest edge in $E(P_1) \sqcup E(P_2)$ (with respect to $<$) is contained in $P_1$.

\begin{definition}[External activity]\label{definition: external activity}
Let $T$ be an arborescence of a strongly connected digraph $G$ rooted at $s \in V$. 
For any edge $e \in E(G) \setminus E(T)$, 
there are exactly two edge-disjoint directed paths $P_1$ and $P_2$
that share  the same starting vertex and ending vertex.
Let $P_1$ be the path that contains $e$.
We say that $e$ is 
\emph{externally active} with respect to $T$ 
if $P_1< P_2$.
The  \emph{external activity} $\ext(T)$  of  $T$ is the number of edges in $G$ that  are externally active with respect to $T$.
\end{definition}

\begin{figure}[ht!]
\begin{tabular}{c c c}
\includegraphics[scale=1]{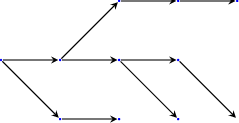} &
\includegraphics[scale=1]{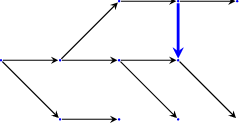} &
\raisebox{0.75\height}{\includegraphics[scale=1]{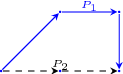}} \\
(a) & (b) & (c)
\end{tabular}
      \caption{(a) An arborescence $T$. (b) An arborescence $T$ with an extra edge $e$.
      (c) The (undashed) path $P_1$ that contains $e$  and the  (dashed) path $P_2$ that doesn'nt contain $e$. }
      \label{figure: external activity}
\end{figure}

See Figure~\ref{figure: external activity} for an illustration describing the process in Definition~\ref{definition: external activity}.

\begin{definition}[Greedoid polynomial]\label{definition: greedoid polynomial}
Let $G$ be a strongly connected digraph, and let $s \in V$.
The  \emph{(single variable) greedoid polynomial}  is 
\begin{equation*}
T(G,s;y):=\sum_{T}  y^{\ext(T)},
\end{equation*}
with the sum  taken over all arborescences of $G$ rooted at $s$.
\end{definition}

This definition of the greedoid polynomial
is due to Bj{\"o}rner,  Korte,  and Lov{\'a}sz~\cite{BKL85}.
Their definition  encompasses a  bigger
 for a more general family of combinatorial objects called \emph{greedoids},
 of  which   the polynomial 
 in Definition~\ref{definition: greedoid polynomial} is a special case.

The polynomial $T(G,s;y)$ does not depend on the choice of total order $<$ on the edges~\cite[Theorem~6.1]{BKL85}. 
If $G$ is a loopless undirected graph,
then the greedoid polynomial $T(G,s;y)$ of $G$
(considered as a bidirected digraph)
 is equal to 
$y^{E(G)}\mathcal T( G;1,y)$,
where $\mathcal T( G;x,y)$ is the Tutte polynomial of $G$ (considered as an undirected graph)~\cite{BKL85}.
We remark that the extra factor $y^{E(G)}$ is due to undirected edges of $G$ being considered as two separate directed edges.
%

We refer the reader  to \cite{BZ92}
for an introduction to greedoids and related topics, and \cite{GM89,GT90,GM91} for a more detailed study of the greedoid polynomial.




\begin{definition}[Reverse $G$-parking functions]\label{definition: reverse G-parking functions}
Let $G$ be a strongly connected digraph, and let $s \in V$.
A \emph{reverse $G$-parking function} with respect to $s$ is a function $f:V\setminus \{s\} \to \N_0$ such that, for any non-empty subset $A \subseteq V \setminus \{s\}$,
there exists $v \in A$ for which $f(v)$ is strictly smaller than the number of edges from $V \setminus A$ to $v$.
\end{definition}
  We use $\Park(G,s)$ to denote the set of  reverse $G$-parking functions rooted at $s$.

 $G$-parking functions were originally defined by Konheim and Weiss \cite{KW66} for complete graphs, and
 were then extended to 
 arbitrary digraphs   by Postnikov and Shapiro \cite{PS04}.
  Reverse $G$-parking functions are known under several different names, including  {reduced divisors}~\cite{BS13},  {superstable configurations}~\cite{HLM08, AB11}, and $\mathbf{\chi}${-superstable configurations}~\cite{GK15}.

We remark that  the choice of working with reverse $G$-parking functions (instead of  $G$-parking functions) is not due to a mere choice of convention,
but is due to a duality property  in the next lemma.

For any sink configuration $\cs$ of $G$,
its \emph{dual function} $f:V \setminus \{s\} \to \N_)$ is given by
\[f(v):= \outdeg(v) -1 -\cs(v) \quad (v \in V\setminus\{s\}). \]

 \begin{lemma}[{\cite[Theorem~4.4]{HLM08}}]\label{l. superstable}
 Let  $G$ be  a connected Eulerian digraph, and let $s \in V$.
  Then a sink configuration of $G$ is sink recurrent if and only if 
  its dual function is a  reverse $G$-parking function. \qed
 \end{lemma}

\begin{remark}
We would like to warn the reader that  Lemma~\ref{l. superstable} 
is false if $G$ is not an Eulerian digraph.
 For an arbitrary (strongly connected) digraph,
the functions dual to recurrent sink configurations are called
  {z-superstable configurations} \cite{AB11, GK15}.
  We refer the reader to \cite[Section~4]{GK15} (specifically, Example~4.17)
  for the subtle distinction between these two functions.
\end{remark}

\begin{definition}[Level of a function]\label{definition: level of a function}
The \emph{level} of a  function $f:V\setminus \{s\} \to \N_0$ is
 \[\lvl(f):= |E|-|V|+1  -\sum_{v \in V \setminus \{s\}} f(v). \qedhere\]
\end{definition}
Note that the level of a function is equal to level of its dual sink configuration plus the outdegree of $s$.

\subsection{Cori-Le Borgne bijection for directed graphs}
In this subsection we give a bijection between reverse $G$-parking functions and arborescences of a directed graph $G$.
This bijection is a directed graph version of  Cori-Le Borgne bijection~\cite{CL03,BS13} for undirected graphs.
For the description of this bijection, see
 Algorithm~\ref{algorithm: parking functions to arborescences}.

\begin{algorithm}
\label{algorithm: parking functions to arborescences}
\caption{Cori-Le Borge bijection from $G$-parking functions to  arborescences of $G$.}
\KwIn{\\$G$-parking function $f$ with respect to $s$, 
\\Total order on the edges of $G$.}
\KwOut{\\Arborescence $T_f$ of $G$ rooted at $s$.}

\BlankLine
{\bf Initialization:}
\\$\BV:=\{s\}$ (burnt vertices), 
\\$\BE:=\emptyset$ (burnt edges), 
\\$T:=\emptyset$ (directed tree).

\BlankLine

 \While{$\BV \neq V(G)$ \label{item: algorithm while start}}{
$e:=\max\{ (v,w) \in E(G) \, | \, (v,w) \not\in \BE,  \, v \in \BV, w \not\in \BV \}$,\\
 $w:=$ the target vertex of  $e$, \\
\If{$f(w) ==$  {the number of edges in $\BE$ with} $w$ as the target vertex \label{item: algorithm parking function comparison} }{
$\BV \leftarrow \BV \cup \{w\}$,
\\ $T \leftarrow T\cup \{e\}$,
}
$\BE \leftarrow \BE \cup \{e\}$ \label{item: algorithm if end}
}
\label{item: algorithm while end}
Output $T_f:=T$.
\end{algorithm}


\begin{reptheorem}{t. generalized Cori-Le Borgne}
Let $G$ be a strongly connected digraph, and let $s \in V$.
Then Algorithm~\ref{algorithm: parking functions to arborescences}
is a bijection that sends   reverse $G$-parking functions with respect to $s$ to arborescences of $G$ rooted at $s$.
Furthermore, the external activity of the output arborescence is the level of the input reverse $G$-parking function.
\end{reptheorem}

\begin{remark}
There are several bijections in the existing literature
between  $G$-parking functions and spanning trees     of undirected graphs (for example \cite{Biggs99t, KY08, Back12, YBKS14,PYY17}).
For directed graphs,
   there is a   
bijection between    reverse $G$-parking functions  and arborescences of $G$ by Chebikin and Pylyavskyy \cite{CP05}.
Note that this bijection is different  from the bijection in Algorithm~\ref{algorithm: parking functions to arborescences} as the former  does not preserve the notion of activities (see \cite[Section~5]{CP05}). 
\end{remark}

The following theorem is a direct consequence of Theorem~\ref{t. generalized Cori-Le Borgne} and Lemma~\ref{l. superstable}.

\begin{reptheorem}{t. Merino's theorem}[Merino's Theorem for Eulerian digraphs]
Let $G$ be a connected Eulerian digraph.
Then for any $s \in V$,
\[
\pushQED{\qed} 
T(G,s;y)=\sum_{\cs \in \Rec(G,s)} y^{\lvl(\cs)+\outdeg(s)}.  \qedhere
\popQED
\]     
 \end{reptheorem}
If the graph in Theorem~\ref{t. Merino's theorem} is bidirected, then we recover the original theorem of Merino L{\'o}pez~\cite{Mer97}.

\begin{remark}
We would like to warn the reader that there are non-Eulerian digraphs for which 
Theorem~\ref{t. Merino's theorem} is false.
This is because $T(G,s;1)$ is the number of  arborescences of $G$,
while $|\Rec(G,s)|$ is the number of reverse arborescences of $G$ (by Lemma~\ref{lemma: sandpile group bijection}\eqref{item: sandpile group 2}).
Those two numbers are in general  not equal for non-Eulerian digraphs.
\end{remark}

The following corollary is a consequence of Theorem~\ref{t. Merino's theorem}.

\begin{corollary}\label{corollary: greedoid polynomial and Biggs-Merino polynomial}
Let $G$ be a connected Eulerian digraph.
Then  for any $s \in V$,
\[ T(G,s;y)=\B(G,s;y). \]
\end{corollary}
\begin{proof}
Since $G$ is a connected Eulerian digraph,
the primitive period vector $\rr$ of $G$ is equal to $(1,\ldots,1)$.
By the markov chain tree theorem~\cite{AT89},
this implies that 
 the period constant $\alpha$ is equal to the number of reverse arborescences rooted at $s$.
By   Proposition~\ref{p. size Rec(G,es)} and Lemma~\ref{lemma: sandpile group bijection}\eqref{item: sandpile group 2}, this implies that there are as many  recurrent sink classes as  recurrent sink configurations.
Hence we conclude that each recurrent sink class of $G$ contains a unique recurrent sink configuration.
Together with Theorem~\ref{t. Merino's theorem}, this implies that
\begin{equation*}
T(G,s;y)=\sum_{\cs \in \Rec(G,s)} y^{\lvl(\cs)+\outdeg(s)}= \sum_{[\cs] \in \Rec(G,\es)} y^{\lvl([\cs])}=\B(G,s;y). \qedhere
\end{equation*}
\end{proof}
This relates the Biggs-Merino polynomial to the greedoid polynomial, as promised in the beginning of this section.

Corollary~\ref{corollary: greedoid polynomial and Biggs-Merino polynomial} gives two interesting consequences for a connected Eulerian digraph $G$.
The first consequence is that  the greedoid polynomial $T(G,s;y)$ does not depend on the choice of $s$ (by Theorem~\ref{t. main theorem}).
The second consequence is that 
 $\B(G,s;2)$ 
counts the number of 
  subgraphs of $G$ such that, for any $v \in V$, there exists a directed path from $s$ to $v$ in the subgraph
   (since $T(G,s;2)$ counts the same thing by \cite[Lemma~2.1]{GM89}).
%
%
%

\begin{remark}
We would like to warn the reader that Corollary~\ref{corollary: greedoid polynomial and Biggs-Merino polynomial} is false when $G$ is a non-Eulerian digraph.
This is because the number of reverse arborescences of $G$ depends on the choice of $s$ if $G$ is non-Eulerian, 
while $\B(G,s;y)$ does not depend on the choice of $s$ (by Theorem~\ref{t. main theorem}).
\end{remark}

The rest of this section is focused on the proof of Theorem~\ref{t. generalized Cori-Le Borgne}.


For any $G$-parking function $f$,
denote by $T_f$ the output of Algorithm~\ref{algorithm: parking functions to arborescences}, 
denote by $\BV(f)$ the set of vertices that are burnt in Algorithm~\ref{algorithm: parking functions to arborescences},
and by $\BE(f)$ the set of edges that are burnt in Algorithm~\ref{algorithm: parking functions to arborescences}.

\begin{lemma}\label{lemma: properties of algorithm 1}
Let $G$ be a strongly connected digraph,  let $s \in V$, and let
$f$ be a $G$-parking function with respect to $s$.
Then:
\begin{enumerate}
\item \label{item: all vertices are burnt} $\BV(f)=V(G)$;  
\item \label{item: burnt edges determine the parking function} $f(v) =|\{e \in \BE(f) \mid  \trgt(e)=v \}|$ for all $v \in V$; and
\item \label{item: output of algorithm 1 is an arborescence} 
$T_f$ is an arborescence of $G$ rooted at $s$.
\end{enumerate}
\end{lemma}
\begin{proof}
\begin{enumerate}
\item 
Suppose to the contrary that Algorithm~\ref{algorithm: parking functions to arborescences} terminates when $\BV(f) \subsetneq V(G)$.
Line \ref{item: algorithm while start}-\ref{item: algorithm while end} of the algorithm imply that all  edges with  source vertex in $\BV(f)$ and target vertex in $V(G) \setminus \BV(f)$ are burnt.
Write $A:=V(G) \setminus \BV(f)$.
Line~\ref{item: algorithm parking function comparison} of the algorithm then implies that for  all $v \in A$, the function $f(v)$
is greater than or equal  to the number of edges $V(G) \setminus A$ to $v$.
This contradicts the assumption that $f$ is a $G$-parking function,
as desired.
\item Since $\BV(f)=V(G)$ by Lemma~\ref{lemma: properties of algorithm 1}\eqref{item: all vertices are burnt},
Line~\ref{item: algorithm parking function comparison} of Algorithm~\ref{algorithm: parking functions to arborescences} implies that 
$f(v)$ is equal to the number of burnt edges with target vertex $v$ for all $v \in V$, as desired.

\item It follows from Line \ref{item: algorithm while start}-\ref{item: algorithm while end}  of Algorithm~\ref{algorithm: parking functions to arborescences} that
$T_f$ is a directed tree with $|\BV(f)|-1$ edges and with $s$ as the unique source vertex.
 Since $\BV(f)=V(G)$ by Lemma~\ref{lemma: properties of algorithm 1}\eqref{item: all vertices are burnt}, it then follows that $T_f$ is an arborescence of $G$ rooted at $s$.
 \qedhere
\end{enumerate}
\end{proof}

\begin{lemma}\label{lemma: external activity becomes level by Cori-Le Borgne bijection}
Let $G$ be a strongly connected digraph, let $s \in V$, and let $f$ be a $G$-parking function with respect to $v$.
Then an edge $e \in E(G) \setminus E(T_f)$ is externally active with respect to $T_f$ if and only if $e$ is not contained in $\BE(f)$.
\end{lemma}

\begin{proof}
Let $P_1$ and $P_2$ be two edge-disjoint directed paths as in Definition~\ref{definition: external activity}.
Note that $e$ is contained in $P_1$ by definition.
Let $e'$ be the minimum edge in $E(P_1) \sqcup E(P_2)$.
We need to show that $e'$ is contained in $P_2$ if and only if 
$e$ is contained in $\BE(f)$.

Suppose that $e'$ is contained in $P_2$.
By the minimality of $e'$,
it then follows that the source vertex of $e$ is burnt before $e'$ in the while loop of Algorithm~\ref{algorithm: parking functions to arborescences}.
Again by the minimality of $e'$,
it then follows that $e$ is evaluated before $e'$ in the while loop of the algorithm.
Since $e$ is not contained in $T_f$, it then follows that $e$ is burnt when it is evaluated.
This proves one direction of the claim.

Suppose that $e'$ is contained in $P_1$.
By the minimality of $e'$, it then follows that all edges in 
$P_2$ are evaluated before $e'$ in the while loop of Algorithm~\ref{algorithm: parking functions to arborescences}.
This implies that all vertices in $P_2$ is burnt before $e'$ is evaluated by the while loop.
Since $P_1$ and $P_2$ share the same target vertex and $e$ is the last edge in $P_1$,
it then follows that $e$ is either not evaluated or evaluated after its target vertex is burnt in the while loop.
In either cases $e$ is not burnt in the while loop.
This proves the other direction of the claim.
\end{proof}

We now give an algorithm that will provide 
 the inverse map to 
 Algorithm~\ref{algorithm: parking functions to arborescences} (note that at   this point we have not yet shown that Algorithm~\ref{algorithm: parking functions to arborescences} is a bijection).
See Algorithm~\ref{algorithm: arborescences to parking functions} for the description of the algorithm.

\begin{algorithm}
\label{algorithm: arborescences to parking functions}
\caption{Cori-Le Borge bijection from  arborescences of $G$ to $G$-parking functions.}
\KwIn{\\Arborescence $T$ of $G$ rooted at $s$, 
\\Total order on the edges of $G$.}
\KwOut{\\$G$-parking function $f_T$ with respect to $s$.}

\BlankLine
{\bf Initialization:}
\\$\BV:=\{s\}$ (burnt vertices), 
\\$\BE:=\emptyset$ (burnt edges).

\BlankLine

 \While{$\BV \neq V(G)$ }{
$e:=\max\{ (v,w) \in E(G) \, | \, (v,w) \not\in \BE,  \, v \in \BV, w \not\in \BV \}$,\\
 $w:=$ the target vertex of  $e$, \\
\If{$e\in E(T)$ }{
$\BV \leftarrow \BV \cup \{w\}$,
}
$\BE \leftarrow \BE \cup \{e\}$ 
}
Output $f_T$, with $f_T(v):=$ the number of edges in $\BE$ with $v$ as target vertex (for $v \in V \setminus \{s\}$).
\end{algorithm}

For any arborescence $T$ of $G$,
denote by $f_T$ the output of Algorithm~\ref{algorithm: arborescences to parking functions}, 
denote by $\BV(T)$ the set of vertices that are burnt in Algorithm~\ref{algorithm: arborescences to parking functions},
and by $\BE(T)$ the set of edges that are burnt in Algorithm~\ref{algorithm: arborescences to parking functions}.

\begin{lemma}\label{lemma: inverse function to Cori-Le Borgne map}
Let $G$ be a strongly connected digraph, let $s \in V$, and let $T$ be an arborescence of $G$ rooted at $s$.
Then:
\begin{enumerate}
\item \label{item: output of arborescences are G-parking functions} $f_T$ is a $G$-parking function with respect to $s$; and
\item \label{item: inverse function Cori-Le Borgne map} For any $G$-parking function $f$ with respect to $s$, we have $f_{T_f}=f$.
\end{enumerate}
\end{lemma}
\begin{proof}
\begin{enumerate}
\item Let $A$ be an arbitrary non-empty subset of $V \setminus \{s\}$.
Since $T$ is an arborescence of $G$, it follows that 
Algorithm~\ref{algorithm: arborescences to parking functions} terminates only when all vertices are burnt.
Let $v$ be the first vertex in $A$ that is burnt by Algorithm~\ref{algorithm: arborescences to parking functions}.
By the minimality assumption on $v$,
it then follows that  the source vertex of every edge in $\{e \in \BE(T) \mid \trgt(e)=v \}$ is contained in $V \setminus A$.
Also note that any edge in $T$ is not burnt in Algorithm~\ref{algorithm: arborescences to parking functions}.
Hence:
\begin{align*}
f_T(v)=&| \{e \in \BE(T) \mid \trgt(e)=v \}| \\
\leq& | \{e \in E(G) \mid  \src(e) \in V\setminus A, \trgt(e)=v, \text{ and } e \notin E(T)\} |\\
=&\, \text{the number of edges from $V \setminus A$ to $v$} \ -1.
\end{align*}
 Since the choice of $A$ is arbitrary, this shows that $f_T$ is a $G$-parking function.

\item It follows from the description of Algorithm~\ref{algorithm: parking functions to arborescences} and Algorithm~\ref{algorithm: arborescences to parking functions} that 
$\BE(f)=\BE(T_f)$.
It then follows from Lemma~\ref{lemma: properties of algorithm 1}\eqref{item: burnt edges determine the parking function} that 
$f=f_{T_f}$. \qedhere
\end{enumerate}
\end{proof}

\begin{proof}[Proof of Theorem~\ref{t. generalized Cori-Le Borgne}]
Note that  Algorithm~\ref{algorithm: parking functions to arborescences} maps $G$-parking functions to arborescences of $G$ (by Lemma~\ref{lemma: properties of algorithm 1}\eqref{item: output of algorithm 1 is an arborescence}).
Also note that Algorithm~\ref{algorithm: arborescences to parking functions}
maps arborescences of $G$ to $G$-parking functions (by Lemma~\ref{lemma: inverse function to Cori-Le Borgne map}\eqref{item: output of arborescences are G-parking functions}).
Finally, note that applying Algorithm~\ref{algorithm: arborescences to parking functions} after Algorithm~\ref{algorithm: parking functions to arborescences} sends a $G$-parking function back to itself.
These three statements imply that 
Algorithm~\ref{algorithm: parking functions to arborescences} is a bijection from $G$-parking functions to arborescences of $G$,
and Algorithm~\ref{algorithm: arborescences to parking functions} is its inverse.

It follows from Lemma~\ref{lemma: properties of algorithm 1}\eqref{item: burnt edges determine the parking function} and Lemma~\ref{lemma: external activity becomes level by Cori-Le Borgne bijection} that 
$\ext(T_f)=\lvl(f)$ for any $G$-parking function $f$.
The proof is complete.
\end{proof}


\section{Concluding remarks}\label{s. conjecture}
In this section we present a few unanswered questions  that might warrant further research.

As remarked in Section \ref{s. greedoid polynomial},
there are non-Eulerian digraphs for which  the conclusion of
Theorem~\ref{t. Merino's theorem} and Corollary~\ref{corollary: greedoid polynomial and Biggs-Merino polynomial} are false.
A further research can be done on extending these two theorems to general digraphs, and here we list two questions of that flavor.

\begin{question}\label{question: greedoid for strongly connected digraphs}
Let $G$ be a non-Eulerian strongly connected digraph.
\begin{enumerate}
\item \label{item: question greedoid} Does there exist 
 a greedoid for which its greedoid polynomial satisfies the conclusion of 
 Theorem~\ref{t. Merino's theorem} or Corollary~\ref{corollary: greedoid polynomial and Biggs-Merino polynomial}?
\item \label{item: question general}  Does there exist an expression for the polynomial in the right hand side of Theorem~\ref{t. Merino's theorem} or Corollary~\ref{corollary: greedoid polynomial and Biggs-Merino polynomial} that is not related to the sandpile model?
\end{enumerate}
\end{question}
Note that Question~\ref{question: greedoid for strongly connected digraphs}\eqref{item: question greedoid} is a special case of Question~\ref{question: greedoid for strongly connected digraphs}\eqref{item: question general}.

One consequence of Merino's Theorem for undirected graphs is that it implies Stanley's pure $O$-sequence conjecture~\cite{Stan96} for cographic matroids.
It is therefore natural to ask for a relationship between our works and $O$-sequences.

Let $X$ be a finite, nonempty set of
 (monic) monomials in the indeterminates $x_1,\ldots, x_k$.
We call $X$ 
  a \emph{(monomial) order ideal} if, for any monomial $m_1 \in X$ and any monomial $m_2$, we have $m_2$ divides $m_1$ only if    $m_2 \in X$.
  An order ideal $X$ is \emph{pure} if all the maximal monomials in $X$ (i.e. polynomials that are not divisible by any other elements in $X$)
  have the same degree.

Let $h_i$  $(i\geq 0)$ denote the number of monomials in $X$ with degree $i$.
The \emph{$h$-vector} of $X$ is the vector 
($h_0,\ldots, h_n$), where $n$ is the maximum degree of monomials in $X$.
An \emph{$O$-sequence} is the $h$-vector of an order ideal, and a \emph{pure $O$-sequence} is the $h$-vector of a pure order ideal.

It follows from Theorem~\ref{t. generalized Cori-Le Borgne} that, for any strongly connected digraph $G$, the nonzero coefficients of its greedoid polynomial, ordered from the highest degree to the lowest degree, is an $O$-sequence.
A further research can be done on extending this observation to other classes of greedoids.

\begin{question}
Do the nonzero coefficients of the greedoid polynomial of a greedoid form an $O$-sequence? If not, what is the class of greedoids for which this property holds?
\end{question}
We remark that there are Eulerian digraphs for which the corresponding $O$-sequence is not pure, for example~\cite[Figure~10]{PP15}.


Another possible research direction is on the method of computing the Biggs-Merino polynomial efficiently.
There is  a variant of the  deletion-contraction  recursion~\cite{BKL85} for the greedoid polynomial and a M\"{o}bius inversion formula~\cite{Perrot-Pham} for the Biggs-Merino polynomial of an Eulerian digraph.
However, 
we are not aware  of any  formulas of the same type for the Biggs-Merino polynomial of non-Eulerian digraphs.
\begin{question}
Does there exist any kind of  deletion-contraction recurrence for the Biggs-Merino polynomial of non-Eulerian digraphs?
\end{question}

\section*{Acknowledgement}
The author would like to thank  Lionel Levine for suggesting the problem to the author and for many helpful discussions,   Farbod Shokrieh for pointing out the connection between the Biggs-Merino polynomial and the greedoid polynomial,
Spencer Backman for introducing a more compact version of Cori-Le Borgne bijection to the author
and for many other comments, and
 Matthew Farrell for helpful comments. 
Last but not the least,  the author would like to thank  anonymous referees for inspiring suggestions and for pointing out the connection between our works and  the $O$-sequence conjecture.
The research was supported by  NSF grant DMS-1243606.

Test.


\bibliographystyle{alpha}
\bibliography{tutte}

\end{document}